\documentclass[letterpaper, 10pt,conference]{article}

\usepackage[utf8]{inputenc}
\usepackage[english]{babel}
\usepackage{graphicx}
\usepackage{epstopdf}
\usepackage{amsmath,amsfonts,amssymb,amsthm}
\usepackage{mathtools}
\usepackage{enumerate}
\usepackage{paralist}

\usepackage{Hfonts,Hcommands,Henvironments}

\usepackage{hyperref}

\hypersetup{
    bookmarks=true,         
    unicode=true,          
    pdftoolbar=true,        
    pdfmenubar=true,        
    pdffitwindow=true,     
    pdfstartview={FitH},    
    pdftitle={Lyapunov functions},    
    pdfauthor={H. Stein Shiromoto},     
    pdfsubject={Subject},   
    pdfcreator={TeXShop},   
    pdfproducer={H. Stein Shiromoto}, 
    pdfkeywords={Lyapunov Functions} {matlab}, 
    pdfnewwindow=true,      
    colorlinks=true,       
    linkcolor=black,          
    citecolor=black,        
    filecolor=black,      
    urlcolor=black           
}

\usepackage[geometry]{ifsym}

\title{\LARGE Interconnecting a System Having a Single Input-to-State Gain With a System Having a Region-Dependent Input-to-State Gain}
\author{Humberto \textsc{Stein Shiromoto}$^{1}$, Vincent \textsc{Andrieu}$^2$, Christophe \textsc{Prieur}$^1$%
\thanks{\footnotesize The authors are with $^1$ Control Systems Department, GIPSA-lab, Grenoble,
France; $^2$ LAGEP, Universit\'e de Lyon 1, France. E-mail: \texttt{
\href{mailto:humberto.shiromoto@ieee.org}{humberto.shiromoto@ieee.org}}. The works of the first and the third authors are partly supported by HYCON2 Network of Excellence Highly-Complex and Networked Control Systems, grant agreement 257462.}}

\begin{document}
\maketitle
\thispagestyle{empty}
\pagestyle{empty}

\begin{abstract}
	For an ISS system, by analyzing local and non-local properties, it is obtained different input-to-state gains. The interconnection of a system having two input-to-state gains with a system having a single ISS gain is analyzed. By employing the Small Gain Theorem for the local (resp. non-local) gains composition, it is concluded about the local (resp. global) stability of the origin (resp. of a compact set). Additionally, if the region of local stability of the origin strictly includes the region attraction of the compact set, then it is shown that the origin is globally asymptotically stable. An example illustrates the approach.
\end{abstract}

\section{Introduction}

	The use of nonlinear input-output gains for the study of the stability of nonlinear was introduced in \cite{Zames1963,Zames1966P1} by considering a system as an input-output operator. The condition that ensures stability, called Small Gain Theorem, of the resulting interconnected system is based on the contraction principle (\cite{Zames1966P1}).

The works \cite{Sontag1989} and \cite{Sontag1990} introduce a new concept of gain relating the input to system states. This notion of stability, called Input-to-State Stability (ISS), combines Zames and Lyapunov approaches (\cite{Sontag:2001,Sontag2008}). Characterizations in terms of dissipation and Lyapunov functions are given in \cite{SontagWang1995} and \cite{SontagWang1996}.

In \cite{Jiangetal1994}, the contraction principle is used in the ISS notion to obtain an equivalent Small Gain Theorem. A formulation of this criteria in terms of Lyapunov functions may be found in \cite{Jiangetal:1996} and \cite{Liberzon:2012}.

Besides stability analysis, the Small Gain Theorem may also be used for the design of dynamic feedback laws satisfying robustness constraints. The interested reader may see \cite{Doyleetal1992,FreemanKokotovic2008,Isidori:1999} and \cite{Sastry:1999} and references therein.

	Other versions of the Small Gain theorem do exist in the literature, examples of which can be found in \cite{AstolfiPraly:2012}. See also \cite{Angeli:2007,Ito2006} and \cite{ItoJiang:2009} for the interconnection of possibly non-ISS systems.

	In order to apply the Small Gain Theorem, it is required that the composition of the nonlinear gains is smaller than the argument for all of its positive values (\cite{Jiangetal:1996,Liberzon:2012}). Such a condition, called Small Gain Condition, restricts the application of the Small Gain Theorem to a composition of well chosen gains.

In this work, it is made use of the Small Gain Theorem in a less conservative way. This new condition ensures the asymptotic stability of a system by showing that if there exist two different gains compositions such that they satisfy the Small Gain Condition, not for all values of the arguments, but in two different regions, and if these regions cover the set of all positive values, then the resulting interconnected system is globally asymptotically stable. Thus, this approach may be seen as a composition of two different small gain conditions that hold in different regions: a local and a global.

The use of a unifying approach is well known in the literature, see \cite{AndrieuPrieur2010} for the combination of control Lyapunov functions and \cite{Chaillet:2012} for a stability concept uniting ISS and the integral variant of ISS (namely, iISS \cite{Sontag1998iISS}) properties.

This paper is organized as follows. In Section \ref{sec:introduction:Preliminaires}, the basic concepts of Input-to-State Stability and the Dini derivative are presented. Also, the system under consideration, the problem statement and a motivational example are presented. In Section \ref{sec:asymmetric:assumptions}, the assumptions to solve the problem under consideration, as well as the main results are presented. Section \ref{sec:asymmetric:illustration} presents an example that illustrates the assumptions and main results. Section \ref{sec:asymmetric:proof} contains the proofs of the main results. Concluding remarks are given in Section \ref{sec:conclusion}. Finally, in Section \ref{sec:auxiliary results}, auxiliary results are stated. Due to space limitations, some of the proofs were omitted.

{\small\textbf{Notation.} Let $\bS$ be a subset of $\R^n$ containing the origin, the notation $\bS_{\neq0}$ stands for $\bS\setminus\{0\}$. The closure of $\bS$ is denoted by $\overline{\bS}$. Let $x\in\R^n$, the notation $|x|$ stands for the Euclidean norm of $x$. A function $f:\bS\to\R$ defined in a subset $\bS$ of $\R^n$ containing $0$ is \emph{positive definite} if, $\forall x\in\bS_{\neq0}$, $f(x)>0$ and $f(0)=0$. It is \emph{proper} if $f(|x|)\to\infty$ as $|x|\to\infty$. By $\cC^k$ it is denoted the class of $k$-times continuously differentiable functions, by $\cK$ it is denoted the class of continuous and strictly increasing functions $\gamma:\Ras\to\Ras$ such that $\gamma(0)=0$; it is denoted by $\cK_\infty$ if, in addition, they are unbounded. Let $c\in\Rs$, the notation $\Omega_{c}(f)$ stands for the subset of $\R^{n}$ defined by $\{x\in\R^n:f(x)<c\}$. Let $x,\bar{x}\in\Ras$, the notation $x\nearrow\bar{x}$ (resp. $x\searrow\bar{x}$) stands for $x\to\bar{x}$ with $x<\bar{x}$ (resp. $x>\bar{x}$).}

\section{Background and problem statement}\label{sec:introduction:Preliminaires}

Consider the system
\begin{equation}\label{eq:x subsystem}
			\dx(t)=f(x(t),u(t)),
\end{equation}
where, $\forall t\in\Ras$, $x(t)\in\R^n$ and $u(t)\in\R^m$, for some positive integers $n$ and $m$. The map $f:\R^n\times\R^m\to\R^n$ is assumed to be continuous, locally Lipschitz on $x$ and uniformly in $u$ on compact sets. A solution of \eqref{eq:x subsystem} with initial condition $x$, and input $u$ at time $t$ is denoted by $X(t,x,u)$. Assume that the origin is an equilibrium point for the system \eqref{eq:x subsystem}, i.e., $f(0,0)=0$.

\begin{definition}
	Consider the function $\xi:[a,b)\to\R$, the limit
	\begin{equation*}
		D^+\xi(t)=\limsup_{\tau\searrow0}\tfrac{\xi(t+\tau)-\xi(t)}{\tau}
	\end{equation*}
	(if it exists) is called \emph{Dini derivative}. Let $k$ be a positive integer. Consider the functions $\varphi:\R^k\to\R$ and $h:\R^k\to\R^k$, the limit
	\begin{equation*}
		D^+_h\varphi(y)=\limsup_{\tau\searrow0}\tfrac{\varphi(y + \tau h(y))-\varphi(y)}{\tau}.
	\end{equation*}
	(if it exists) is called \emph{Dini derivative of $\varphi$ in the $h$-direction at $y$}.
\end{definition}

\begin{definition}\label{def:ISS-Lyapunov function}
	A continuous locally Lipschitz function $V:\R^n\to\R$ is called an ISS-Lyapunov function for system \eqref{eq:x subsystem} if
	
 \noindent$\bullet$ there exist class $\cK_\infty$ functions $\underline{\alpha}$ and $\overline{\alpha}$ such that, $\forall x\in\R^n$, $\underline{\alpha}(|x|)\leq V(x)\leq\overline{\alpha}(|x|)$;
	
 \noindent$\bullet$ there exist a class $\cK$ function $\alpha_x$, called \emph{ISS gain}, and a continuous positive definite function $\lambda_x:\R^n\to\R$ such that, $\forall (x,u)\in\R^n\times\R^m$,
		\begin{equation}\label{eq:ISS Lyapunov inequality}
			|x|\geq\alpha_x(|u|)\Rightarrow D_f^+V(x,u)\leq-\lambda_x(x)
		\end{equation}
		holds.
\end{definition}

	
	From now on, $V$ will be assumed to be an ISS-Lyapunov function for \eqref{eq:x subsystem}.


Consider the system
\begin{equation}\label{eq:z subsystem}
			\dot{z}(t)=g(v(t),z(t)),
\end{equation}
where, $\forall t\in\Ras$, $v(t)\in\R^n$ and $z(t)\in\R^m$, for some positive integers $n$ and $m$. The map $g:\R^n\times\R^m\to\R^n$ is assumed to be continuous, locally Lipschitz on $z$ and uniformly in $v$ on compact sets. A solution of \eqref{eq:z subsystem} with initial condition $z$, and input $v$ at time $t$ is denoted by $Z(t,z,v)$. Assume that the origin is an equilibrium point for the system \eqref{eq:z subsystem}, i.e., $g(0,0)=0$. Consider also the following
\begin{assumption}\label{hyp:asymmetric:z:ISS}
	There exists a continuous locally Lipschitz function $W:\R^m\to\R$ that is an ISS-Lyapunov function for the $z$-subsystem. More precisely, there exist class $\cK_\infty$ functions $\underline{\beta}$ and $\overline{\beta}$ satisfying, $\forall z\in\R^m$, $\underline{\beta}(|z|)\leq W(z)\leq\overline{\beta}(|z|)$. Furthermore, there exist a class $\cK$ function $\delta$ and a continuous positive definite function $\lambda_z:\R^m\to\R$ such that, $\forall(x,z)\in\R^{n}\times\R^{m}$,
\begin{equation}\label{eq:asymmetric:D_g^+W}
W(z)\geq\delta(V(x))\Rightarrow D_g^+W(x,z)\leq-\lambda_z(z),
\end{equation}
where $V$ is the ISS-Lyapunov function of $x$-subsystem.
\end{assumption}

{\bfseries System under consideration.} Interconnecting systems \eqref{eq:x subsystem} and \eqref{eq:z subsystem} by linking the state of \eqref{eq:x subsystem} with the input of \eqref{eq:z subsystem} and vice versa leads to the following system
\begin{equation}\label{eq:general system}
	\left\{\begin{array}{rcl}
			\dx&=&f(x,z)\\
			\dz&=&g(x,z).
	\end{array}\right.
\end{equation}
Since $f(0,0)=0$ and $g(0,0)=0$, the origin is an equilibrium point for \eqref{eq:general system}. Considering the ISS-Lyapunov inequalities, after the interconnection the following implications
\begin{equation*}
	\begin{array}{rcrcrcl}
		V(x)&\geq&\gamma(W(z))&\Rightarrow&D_f^+V(x,z)&\leq&-\lambda_x(x),\\
		W(z)&\geq&\delta(V(x))&\Rightarrow&D_g^+W(x,z)&\leq&-\lambda_z(z)
	\end{array}
\end{equation*} 
are obtained with suitable class $\cK$ functions $\gamma$ and $\delta$.

A sufficient condition that ensures stability of \eqref{eq:general system} is given by the following
\begin{theorem}\label{thm:SGT}\cite{Jiangetal:1996}
	If,
\begin{equation}\label{eq:introduction:general sgc}
	\forall s\in\Rs,\quad\gamma\circ\delta(s)<s.
\end{equation}
Then, the origin is globally asymptotically stable for \eqref{eq:general system}.
\end{theorem}

{\bfseries Problem statement.} At this point, it is possible to explain the problems that are dealt with, in this work.

 \noindent$\bullet$ {\itshape ISS gains computation}. Although the use of ISS gains renders the analysis of stability easy to work with, it is not a trivial task to compute those gains;
 	
 \noindent$\bullet$ {\itshape Small gain condition}. Since the ISS gain is not unique, it might not be an easy task to find two ISS gains: one for the $x$-subsystem of \eqref{eq:general system} and another for the $z$-subsystem of \eqref{eq:general system} such that their composition satisfies \eqref{eq:introduction:general sgc}, for all positive values of the argument.
An illustration of the problem that is dealt with is presented in the following

\begin{figure}[htpb!]
		\centering
		\includegraphics[width=\textwidth]{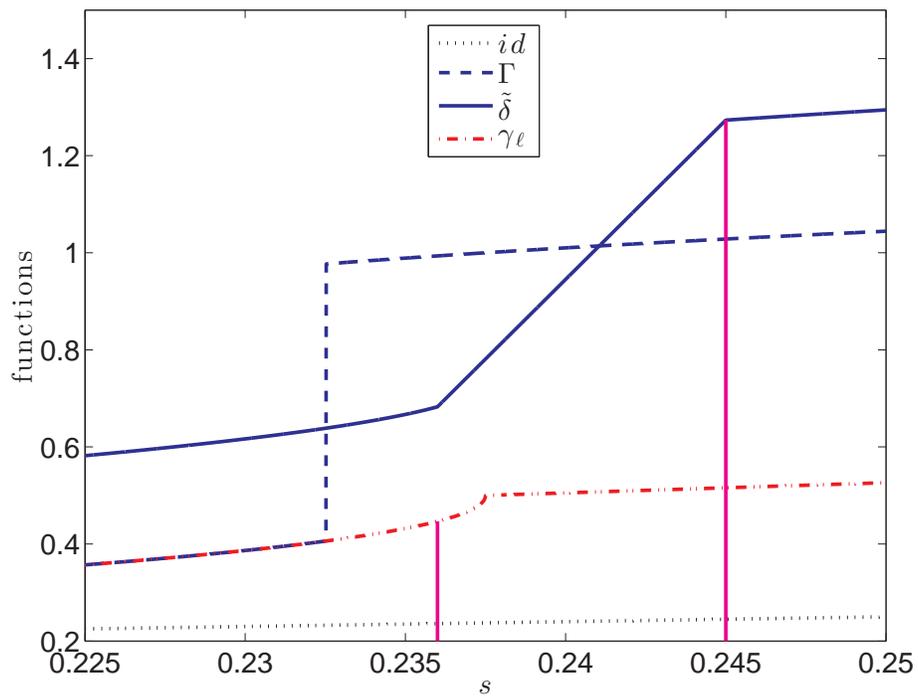}
		\caption{Plot of the functions  $id$ (dotted black line), $\Gamma$ (dashed blue line), $\gamma_\ell$ (dash dotted red line) and the continuous function $\tilde{\delta}$ (solid blue line), in the interval $[0.225,0.25]$. The vertical lines are the values $M_\ell=0.236$ and $M_g=0.245$, respectively.}
		\label{fig:gains}
	\end{figure}

\begin{example}\label{example:asymmetric:introduction}
	Let the functions $f,g:\R\times\R\to\R$ and consider the system 
	\begin{equation}\label{eq:example:general system}
		\left\{\begin{array}{rcccl}
			\dx&=&f(x,z)&=&-\rho(x)+z\\
			\dot{z}&=&g(x,z)&=&-\sign(z)\tilde{\delta}(|z|)+x,
		\end{array}\right.
	\end{equation}
	where $\tilde{\delta}$ will be defined below. Let, $\forall x\in\R$, $V(x)=|x|$, $\rho(x)= 5\pfrac{x}{4}-2x^2+x^3$ and, $\forall z\in\R$, $W(z)=|z|$.
	
	Taking the Dini derivative of $V$ in the $f$-direction, $\forall (x,z)\in\R\times\R$, it yields
	\begin{equation}\label{eq:example:D_f^+V}
			D_f^+V(x,z)\leq-\rho(V(x))+W(z).
	\end{equation}
	This implies that $\exists\varepsilon_x\in(0,1)$ such that, $\forall (x,z)\in\R\times\R$,
	\begin{equation}\label{eq:example: rho(V)>W -> D_f^+V < 0}
		\rho(V(x))\geq \tfrac{W(z)}{1-\varepsilon_x}\Rightarrow D_f^+V(x,z)\leq -\lambda_x(x),
	\end{equation}
	where $\lambda_x(\cdot):=\varepsilon_x\rho(V(\cdot))$. From now on, let $\varepsilon_x=0.05$. Note also that, in the interval $[\pfrac{1}{2},\pfrac{5}{6}]$, $\rho$ is decreasing.

	Consider the piecewise continuous function $\Gamma$ defined by
	\begin{equation}\label{eq:example:Gamma definition}
		\Gamma(s)=\left\{\begin{array}{rcl}
							\rho^{-1}\left(\tfrac{s}{0.95}\right),&\text{if}&s\in\left[0,0.95\rho\left(\tfrac{5}{6}\right)\right),\\
							\rho^{-1}_+\left(\tfrac{s}{0.95}\right),&\text{if}&s\in\left[0.95\rho\left(\tfrac{5}{6}\right),\infty\right),
						 \end{array}
				  \right.
	\end{equation}
	where $[\pfrac{5}{6},\infty)\ni s\mapsto\rho_+(s)=\rho(s)\in[\rho(\pfrac{5}{6}),\infty)$.
	
	\begin{remark}\label{claim:V>Gamma -> D_f^+V<0} The function $\Gamma$ can be viewed as a discontinuous input-to-state gain of the $x$-subsystem of \eqref{eq:example:general system}. More preciselly, $\forall (x,z)\in\R\times\R$, $V(x)\geq\Gamma(W(z))\Rightarrow D_f^+V(x,z)\leq-\lambda_x(x)$. Furthermore, the function $\Gamma$ is ``optimal'', in the sense that if there exist a function $\Gamma^\ast:\R\to\R$ and a value $s^\ast\in\Rs$ such that $\Gamma^\ast(s^\ast)<\Gamma(s^\ast)$, then $\exists(x^\ast,z^\ast)\neq(0,0)$ such that $V(x^\ast)\geq\Gamma^\ast(W(z^\ast))$ and $D_f^+V(x^\ast,z^\ast)>0$.
	\end{remark}
	
	It follows from Remark \ref{claim:V>Gamma -> D_f^+V<0} that an ISS gain for the $x$-subsystem of \eqref{eq:example:general system} is any class $\cK$ function $\gamma$ such that, $\forall s\in\Rs$, $\Gamma(s)\leq \gamma(s)$.
	
	\noindent{\itshape A local gain.} Consider the function $[0,\pfrac{1}{2})\ni s\mapsto \rho_-(s)=\rho(s)\in[0,\rho(\pfrac{1}{2}))$. Since $\rho_-$ is strictly increasing on its domain, it is invertible. Let $\gamma_\ell$ be a class $\cK$ function such that, $\forall s\in[0,0.95\rho(\pfrac{1}{2}))$,
	\begin{equation}\label{eq:example:definition of gamma_l}
		\gamma_\ell(s)=\rho_-^{-1}\left(\tfrac{s}{0.95}\right).
	\end{equation}
		Note that, $\gamma_\ell$ satisfies the following inequalities
		\begin{equation*}
			\begin{array}{rl}
				\forall s\in\left[0,0.95\rho\left(\tfrac{1}{2}\right)\right),&\gamma_\ell(s)\leq\Gamma(s),\\
				\forall s\in\left(0.95\rho\left(\tfrac{5}{6}\right),0.95\rho\left(\tfrac{1}{2}\right)\right),&\gamma_\ell(s)<\Gamma(s).
			\end{array}
		\end{equation*}
	Moreover, $\forall (x,z)\in\Omega_{\rho\left(1/2\right)}(V)\times\R$, 
		\begin{equation}\label{claim:V>gamma_l -> D_f^+V<0}
			V(x)\geq\gamma_\ell(W(z))\Rightarrow D_f^+V(x,z)\leq -\lambda_x(x).
		\end{equation}
	
		Let the constant values $M_\ell=0.236$ and $M_g=0.245$. At this point, it is possible to define the function $\tilde{\delta}$ of the $z$-subsystem of \eqref{eq:example:general system}. It is a function of class $\cK_\infty$ satisfying the following inequalities
		\begin{eqnarray}
			\forall s\in(0,M_\ell],&&\gamma_\ell(s)<\tilde{\delta}(s),\label{eq:example:delta gamma_l < s}\\
			\forall s\in[M_g,\infty),&&\Gamma(s)<\tilde{\delta}(s),\label{eq:example:delta Gamma < s}\\
			\forall s\in\left(\rho\left(\tfrac{5}{6}\right),M_\ell\right),&&\tilde{\delta}(s)<\Gamma(s).\label{eq:example:tilde delta < Gamma}
		\end{eqnarray}
		Equations \eqref{eq:example:delta gamma_l < s} and \eqref{eq:example:delta Gamma < s} correspond to two different small gain conditions, the first may be seen as a small gain condition for small values of the argument while the last as a small gain condition for large values of the argument. Note that \eqref{eq:example:tilde delta < Gamma} implies that Theorem \ref{thm:SGT} cannot be applied.\footnote{To see this fact, note that $M_\ell<0.95\rho(\pfrac{1}{2})$. Since $\tilde{\delta}$ is of class $\cK_\infty$ and from \eqref{eq:example:tilde delta < Gamma}, $\forall s\in (0.95\rho(\pfrac{5}{6}),M_\ell)$, $s<\tilde{\delta}^{-1}\circ\Gamma(s)$. Thus, there exists no class $\cK$ function $\gamma$ such that \eqref{eq:introduction:general sgc} holds.} 
		
		Fig. \ref{fig:gains} shows a plot of the functions $\rho$, $\id$, $\Gamma$, $\gamma_\ell$ and $\tilde{\delta}$. 		\end{example}
		
In this work, it will be shown that, if

 \noindent$\bullet$ there exist two ISS gains $\gamma_\ell$ and $\gamma_g$, for the $x$-subsystem of \eqref{eq:general system};
	
 \noindent$\bullet$ there exists one ISS gain $\delta$, for the $z$-subsystem of \eqref{eq:general system};
	
 \noindent$\bullet$ the compositions $\gamma_\ell\circ\delta$ and $\gamma_g\circ\delta$ satisfy the Small Gain Condition, not for all values of the arguments, but for two different intervals ($\bI_\ell,\bI_g\subset\Ras$). In other words,
	\begin{equation*}
			\forall s\in\bI_\ell\setminus\{0\},\gamma_\ell\circ\delta(s)<s\ \text{and},\ 
			\forall s\in\bI_g\setminus\{0\},\gamma_g\circ\delta(s)<s;
	\end{equation*}
	
 \noindent\hspace{2pt}$\bullet$ these intervals are such that $\bI_\ell\cap\bI_g\neq\emptyset$ and $\bI_\ell\cup\bI_g=\Ras$;\\
then, the origin is globally asymptotically stable for \eqref{eq:general system}. See Theorem \ref{prop:asymetric:0 globally stable} below for a precise statement of this result.

\section{Assumptions and main results}\label{sec:asymmetric:assumptions}

In this section, it is specified the assumptions on the system \eqref{eq:general system} necessary to solve the problem under consideration. The proof of the stabilization results are provided from Section \ref{sec:Proof of Proposition 1} to Section \ref{sec:proof of Theorem 2}.

\subsection{Local set of assumptions on the $x$-subsystem}

In this section, it is introduced the set of assumptions to ensure that the origin is locally asymptotically stable for \eqref{eq:general system}.

\begin{assumption}\label{hyp:asymmetric:Lyapunov}
There exist a class $\cK$ function $\gamma_\ell$ and a strictly positive constant $M_{\ell}$ such that, 
\begin{equation}\label{eq:asymmetric:b_l}
		M_\ell<\lim_{s\to\infty}\gamma_{\ell}(s)=b_\ell.
	\end{equation}
Moreover, $\forall (x,z)\in\overline{\Omega_{M_{\ell}}(V)}\times\R^m$,
	\begin{equation}\label{eq:asymmetric:D_f^+V<0:local}
		V(x)\geq\gamma_\ell(W(z))\Rightarrow D_f^+V(x,z)\leq-\lambda_x(x).
	\end{equation}
\end{assumption}

\begin{assumption}\label{hyp:asymmetric:small gain assumption:local}
	The composition of the functions $\gamma_{\ell}$ and $\delta$ is such that,
\begin{equation}\label{eq:asymmetric:local sgc}
\forall s\in(0,M_{\ell}],\quad\gamma_{\ell}\circ\delta(s)<s.
\end{equation}
\end{assumption}

\begin{proposition}\label{prop:asymetric:0 locally stable}
	Under Assumptions \ref{hyp:asymmetric:z:ISS}, \ref{hyp:asymmetric:Lyapunov} and \ref{hyp:asymmetric:small gain assumption:local} the origin is locally asymptotically stable for system \eqref{eq:general system}.
\end{proposition}

\subsection{Non-local set of assumptions on the $x$-subsystem}

In this section, it is introduced the set of assumptions to ensure that a neighborhood of the origin is globally attractive for \eqref{eq:general system}.

\begin{assumption}\label{hyp:asymmetric:Lyapunov global}
	There exist a class $\cK$ function $\gamma_{g}$ and a strictly positive constant $M_{g}$ such that 
	\begin{equation}\label{eq:asymmetric:b_g}
		M_g<\lim_{s\to\infty}\gamma_{g}(s)=b_g.
	\end{equation}
	Moreover, $\forall (x,z)\in(\R^{n}\setminus\Omega_{M_{g}}(V))\times\R^{m}$,
		\begin{equation}\label{eq:asymmetric:D_f^+V<0:global}
		V(x)\geq\gamma_g(W(z))\Rightarrow D_f^+V(x,z)\leq-\lambda_x(x).
	\end{equation}
\end{assumption}

\begin{assumption}\label{hyp:asymmetric:small gain assumption:global}
	The composition of the functions $\gamma_{g}$ and $\delta$ is such that,
\begin{equation}\label{eq:asymmetric:global sgc}
\forall s\in[M_{g},\infty),\quad\gamma_{g}\circ\delta(s)<s.
\end{equation}
\end{assumption}

\begin{proposition}\label{prop:asymetric:neighborhood globally attractive}
	Under Assumptions \ref{hyp:asymmetric:z:ISS}, \ref{hyp:asymmetric:Lyapunov global} and \ref{hyp:asymmetric:small gain assumption:global}, there exist a proper definite positive function $U_g$ and a positive constant $\tilde{M}_g$ such that the set $\overline{\Omega_{\tilde{M}_g}(U_g)}$ is globally asymptotically stable for system \eqref{eq:general system}.
\end{proposition}

\subsection{Main result}

In this section, it is introduced the assumption to ensure that the origin is globally asymptotically stable for \eqref{eq:general system}.

\begin{assumption}\label{hyp:asymmetric:M_l<M_g}
	The positive constants $M_\ell$ and $M_g$ given, respectively, by Assumptions  \ref{hyp:asymmetric:Lyapunov} and \ref{hyp:asymmetric:Lyapunov global} satisfy $M_g<M_\ell$.
\end{assumption}

\begin{theorem}\label{prop:asymetric:0 globally stable}
	Under Assumptions \ref{hyp:asymmetric:z:ISS}-\ref{hyp:asymmetric:M_l<M_g}, the origin is globally asymptotically stable for system \eqref{eq:general system}.
\end{theorem}


\section{Illustration}\label{sec:asymmetric:illustration}

\begin{figure}[htpb!]
		\centering
		\includegraphics[width=0.45\linewidth]{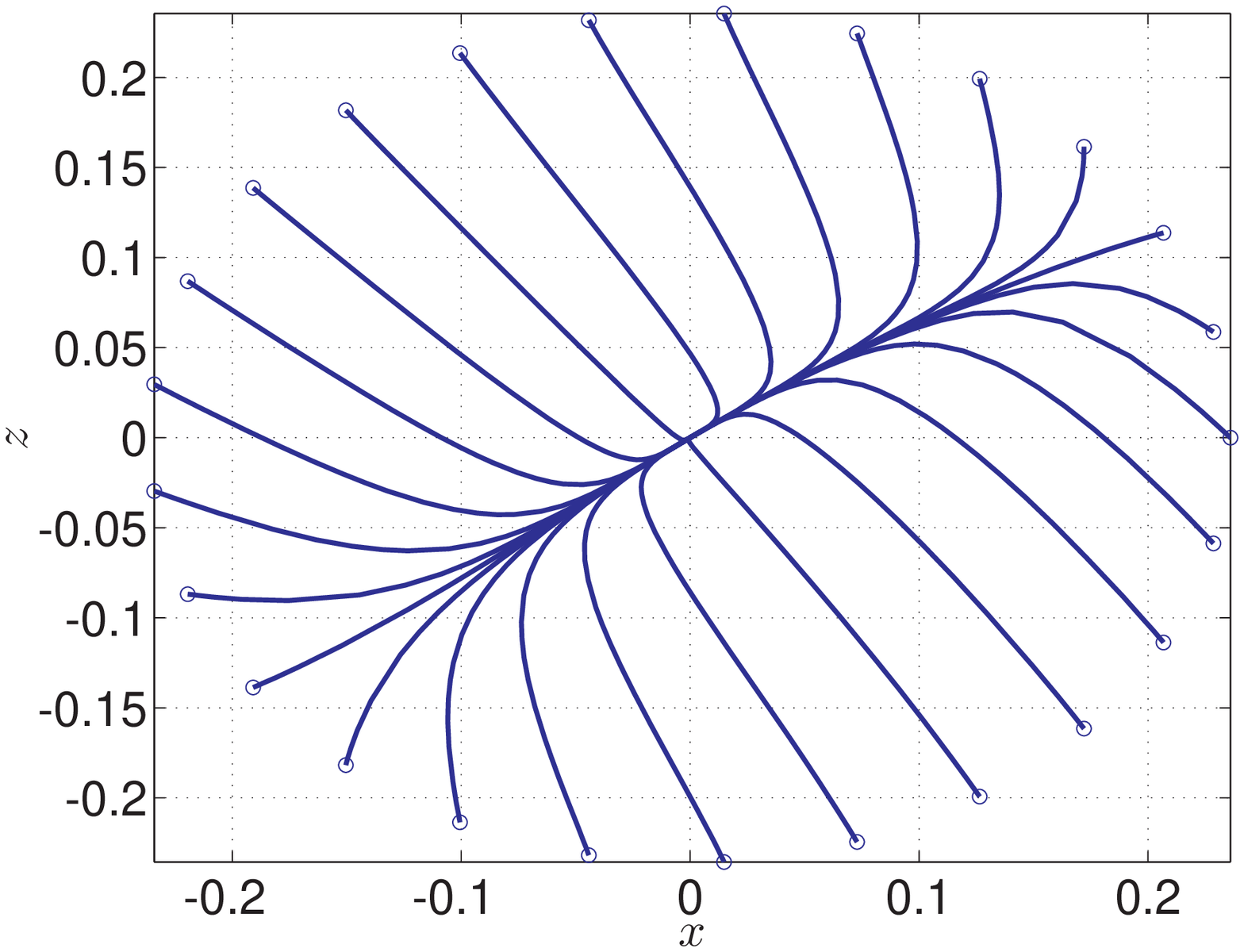}							\includegraphics[width=0.45\linewidth]{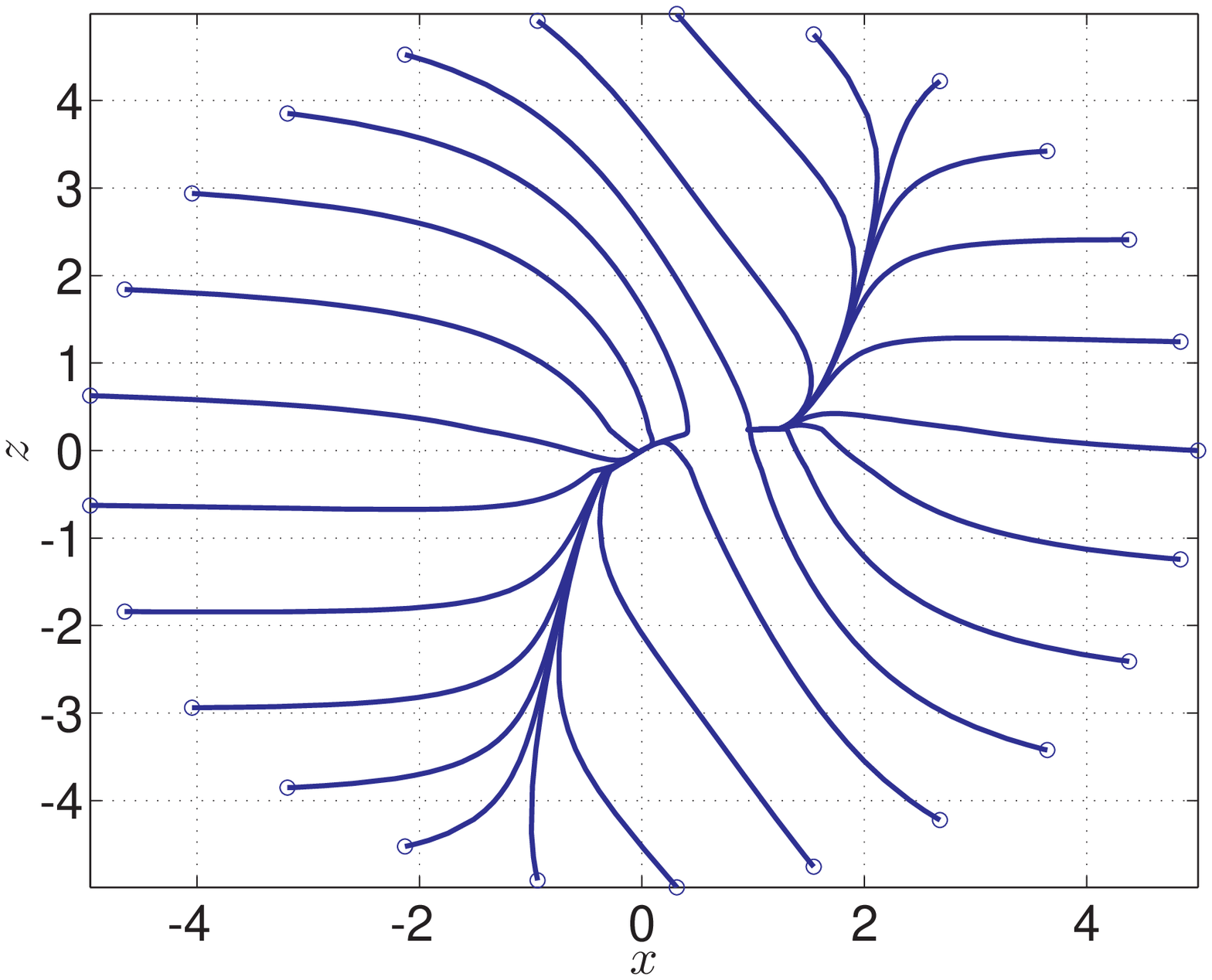}
		\caption{Solutions to \eqref{eq:example:general system} for initial condition starting at a ball centered at the origin with radius, respectively, given by $M_\ell$ and 5.}
		\label{sim:local}
	\end{figure}

\begin{example} \ [Example \ref{example:asymmetric:introduction} revisited.]

	\noindent{\itshape Verifying Assumption \ref{hyp:asymmetric:z:ISS}}. Let the function $\delta$ be given by the inverse of $\tilde{\delta}$. It follows that, $\forall (x,z)\in\R\times\R$, $W(z)\geq \delta(V(x))\Rightarrow D^+_gW(z)\leq-\lambda_z(z)$, where for a given $\varepsilon_z\in(0,1)$ and $\forall z\in\R$, $\lambda_z(z)=\varepsilon_z W(z)$. Thus, Assumption \ref{hyp:asymmetric:z:ISS} holds.
	
\noindent{\itshape Verifying Assumption \ref{hyp:asymmetric:Lyapunov}}. The function $\gamma_\ell$ is given by \eqref{eq:example:definition of gamma_l} and $M_\ell=0.236$. Moreover, it follows from \eqref{claim:V>gamma_l -> D_f^+V<0} that Assumption \ref{hyp:asymmetric:Lyapunov} holds.
	
	\noindent {\itshape Verifying Assumption \ref{hyp:asymmetric:small gain assumption:local}.} It follows from inequality \eqref{eq:example:delta gamma_l < s} that Assumption \ref{hyp:asymmetric:small gain assumption:local} holds. 
	
	From Proposition \ref{prop:asymetric:0 locally stable}, it follows that the origin is locally asymptotically stable for \eqref{eq:example:general system}. Figure \ref{sim:local} shows some solutions of \eqref{eq:example:general system}.
\end{example}

\begin{example} 
	\ [Example \ref{example:asymmetric:introduction} revisited.]
	
	\noindent{\itshape Verifying Assumption \ref{hyp:asymmetric:Lyapunov global}.} Let a class $\cK_\infty$ function $\gamma_g$ be such that, $\forall s\in[0.95\rho(5/6),\infty)$, $\gamma_g(s)=\Gamma(s)$.
	Moreover, $M_g=0.245$. It follows from Remark \ref{claim:V>Gamma -> D_f^+V<0} that, $\forall (x,z)\in(\R\setminus\Omega_{M_g}(V))\times\R$, $V(x)\geq\gamma_g\left(W(z)\right)\Rightarrow D_f^+V(x,z)\leq-\lambda_x(x)$. Thus, Assumption \ref{hyp:asymmetric:Lyapunov global} holds.

{\itshape Verifying Assumption \ref{hyp:asymmetric:small gain assumption:global}.} It follows from inequality \eqref{eq:example:delta Gamma < s} that Assumption \ref{hyp:asymmetric:small gain assumption:global} holds. 

From Proposition \ref{prop:asymetric:neighborhood globally attractive}, it follows that a neighborhood of the origin is globally asymptotically stable for \eqref{eq:example:general system}. Figure \ref{sim:local} shows some solutions of \eqref{eq:example:general system}.
\end{example}


\section{Proofs}\label{sec:asymmetric:proof}

\subsection{Proof of Proposition \ref{prop:asymetric:0 locally stable}}\label{sec:Proof of Proposition 1}

\begin{proof}
	This proof is divided into three parts. In the first one, it is obtained a function $\sigma_\ell$ that is class $\cK_\infty$ and $\cC^1$ with strictly positive derivative. This function is used in the second part, where a class $\cC^0$ proper positive definite function $U_\ell$ is defined and its Dini derivative is studied. In the last part, it is shown that $U_\ell$ is locally Lipschitz and the local asymptotical stability of the origin is concluded by using Lemma \ref{thm:Dini derivative:stability by Dini derivative}.

	{\itshape First part.} Consider the class $\cK$ functions $\delta$ and $\gamma_\ell$ from Assumptions \ref{hyp:asymmetric:z:ISS} and \ref{hyp:asymmetric:Lyapunov}. Under Assumption \ref{hyp:asymmetric:small gain assumption:local}, $\delta$ and $\gamma_\ell$ are such that, $\forall s\in(0,M_\ell]$, $\delta(s)<\gamma_{\ell}^{-1}(s)$.
	
	Since $\gamma_\ell$ is class of $\cK$, from Lemma \ref{lem:construction of K_inf function}, there exists a class $\cK_\infty$ function $\tilde{\gamma}_\ell$ such that, $\forall s\in\Rs$,
	\begin{equation}\label{eq:proof:prop 0 locally stable:delta < t gamma}
		\delta(s)<\tilde{\gamma}_\ell(s)
	\end{equation}
	and, $\forall s\in(0,M_\ell]$,
	\begin{equation}\label{eq:proof:prop 0 locally stable:t gamma < gamma^-1}
		\tilde{\gamma}_\ell(s)<\gamma_\ell^{-1}(s).
	\end{equation}
		
	Since $\delta$ is of class $\cK$ and $\tilde{\gamma}_\ell$ is of class $\cK_\infty$ satisfying, $\forall s\in\Rs$, inequality \eqref{eq:proof:prop 0 locally stable:delta < t gamma}, from Lemma \cite[Lemma A.1]{Jiangetal:1996}, there exists a class $\cK_\infty$ and $\cC^1$ function $\sigma_\ell$ whose derivative is strictly positive and satisfies, $\forall s\in\Rs$,
\begin{equation}\label{eq:proof:prop 0 locally stable:delta < sigma_l < gamma^-1}
	\delta(s)<\sigma_\ell(s)<\tilde{\gamma}_\ell(s).
\end{equation}

	{\itshape Second part.} Let, $\forall(x,z)\in\R^n\times\R^m$, $U_\ell(x,z)\coloneqq\max\{\sigma_\ell(V(x)),W(z)\}$. Note that the function $U_\ell$ is proper positive definite. Pick $(x,z)\in\R^n\times\R^m$, one of three cases is possible: $\sigma_\ell(V({x}))<W({z})$, $W({z})<\sigma_\ell(V({x}))$ or $W({z})=\sigma_\ell(V({x}))$. The proof follows by showing that the Dini derivative of $U_\ell$ is negative. For each case, assume that $(x,z)\in\overline{\Omega_{M_\ell}(V)}\times\R^m$.
	
	\hspace{2pt}\underline{Case 1}. Suppose that $\sigma_\ell(V(x))<W(z)$.
	This implies that $U_\ell(x,z)=W(z)$ and $D^+_{f,g}U_\ell(x,z)=D_g^+W(x,z)$.

	From \eqref{eq:proof:prop 0 locally stable:delta < sigma_l < gamma^-1}, the following inequality $\delta(V(x))<\sigma_\ell(V(x))<W(z)$ holds. Together with \eqref{eq:asymmetric:D_g^+W}, it follows that $D_g^+W(x,z)\leq-\lambda_z(z)$. This concludes Case 1.
	
	\hspace{2pt}\underline{Case 2}. Suppose that $W(z)<\sigma_\ell(V(x))$. 
	This implies that
	$U_\ell(x,z)=\sigma_\ell(V(x))$ and $D^+_{f,g}U_\ell(x,z)=\sigma_\ell'(V(x))D_f^+V(x,z)$.
	From \eqref{eq:proof:prop 0 locally stable:delta < sigma_l < gamma^-1}, the following inequality $W(z)<\sigma_\ell(V(x))<\tilde{\gamma}_\ell(V(x))$ holds. Since $V(x)\leq M_\ell$, it follows that
			\begin{equation}\label{eq:proof:prop 0 locally stable:W<delta<t gamma_l<gamma_l^-1}
					W(z)<\sigma_\ell(V(x))<\tilde{\gamma}_\ell(V(x))<\gamma^{-1}_\ell(V(x)),
			\end{equation}
			where the last inequality follows from \eqref{eq:proof:prop 0 locally stable:t gamma < gamma^-1}. Equation \eqref{eq:asymmetric:D_f^+V<0:local} together with \eqref{eq:proof:prop 0 locally stable:W<delta<t gamma_l<gamma_l^-1} yields $D_f^+V(x,z)\leq-\lambda_x(x)$.
		
		Since, $\forall s\in\Rs$, $\sigma_\ell'(s)>0$, it follows that $D^+_{f,g}U_\ell(x,z)=\sigma_\ell'(V(x))D_f^+V(x,z)\leq-\sigma_\ell'(V(x))\lambda_x(x)$. This concludes Case 2.
		
		\hspace{2pt}\underline{Case 3}. Let $W(z)=\sigma_\ell(V(x))\coloneqq U_\ell^\ast(x,z)$.
		This implies
			\begin{equation*}
				\begin{split}
					D^{+}_{f,g}U_{\ell}^{\ast}(x,z)=\displaystyle\limsup_{t\searrow0}\tfrac{1}{t}(\max\{\sigma_{\ell}(V(X(x,z,t))),\\
					W(Z(z,x,t))\}-U_{\ell}^{\ast}(x,z))\\
					=\displaystyle\limsup_{t\searrow0}\max\left\{\tfrac{\sigma_{\ell}(V(X(x,z,t)))-\sigma_{\ell}(V(x))}{t},\tfrac{W(Z(z,x,t))-W(z)}{t}\right\}\\
					=\max\{\sigma_{\ell}'(V(x))D_f^+V(x,z),D_g^+W(x,z)\}.
									\end{split}
			\end{equation*}
			
			The analysis of $D^+_{f,g}U_\ell^\ast$ is divided in two sub cases. In the first one, the function $D_g^+W$ is analyzed while in the last, the function $D_f^+V$ is analyzed.
			
			\hspace{2pt}\underline{Case 3.a}. {\itshape The analysis of $D_g^+W$.} From \eqref{eq:proof:prop 0 locally stable:delta < sigma_l < gamma^-1}, the following inequality $\delta(V(x))<\sigma_\ell(V(x))=W(z)$ holds. Together with Equation \eqref{eq:asymmetric:D_g^+W}, it yields $D_g^+W(x,z)\leq-\lambda_z(z)$.

			\hspace{2pt}\underline{Case 3.b}. {\itshape The analysis of $D_f^+V$.} From \eqref{eq:proof:prop 0 locally stable:delta < sigma_l < gamma^-1}, the following inequality $W(z)=\sigma_\ell(V(x))<\tilde{\gamma}_\ell(V(x))$ holds. Since $V(x)\leq M_\ell$, it follows that 
			\begin{equation}\label{eq:proof:prop 0 locally stable:W=sigma<t gamma_l<gamma_l^-1}
					W(z)=\sigma_\ell(V(x))<\tilde{\gamma}_\ell(V(x))<\gamma^{-1}_\ell(V(x)),
			\end{equation}
			where the last inequality is due to \eqref{eq:proof:prop 0 locally stable:t gamma < gamma^-1}.
			
		Equation \eqref{eq:asymmetric:D_f^+V<0:local} together with \eqref{eq:proof:prop 0 locally stable:W=sigma<t gamma_l<gamma_l^-1} yields $D_f^+V(x,z)\leq-\lambda_x(x)$.

			To conclude Case 3, $W(z)=\sigma_\ell(V(x))\Rightarrow D^{+}_{f,g}U_{\ell}^{\ast}(x,z)\leq-\min\{\sigma_\ell'(V(x))\lambda_x(x),\lambda_z(z)\}$ holds, since $(x,z)\in\overline{\Omega_{M_\ell}(V)}\times\R^m$.

Let $\tilde{M}_\ell\coloneqq\max\{c\in\Rs:\overline{\Omega_c(U_\ell)}\subset \overline{\Omega_{M_\ell}(V)}\times\{0\}\ \text{and}\ \overline{\Omega_c(U_\ell)}\ \text{is connected}\}$. To sum up all the above cases, $\forall (x,z)\in\overline{\Omega_{\tilde{M}_\ell}(U_\ell)}$,
	\begin{equation}\label{eq:proof:prop 0 locally stable:D+<-max}
		U_\ell(x,z)\leq\tilde{M}_\ell\Rightarrow D^{+}_{f,g}U_{\ell}(x,z)\leq -E_\ell(x,z),
	\end{equation}
	where $E(\cdot,\cdot):=\min\{\sigma_\ell'(V(\cdot))\lambda_x(\cdot),\lambda_z(\cdot)\}$  is continuous and positive definite.
	
	{\itshape Third part.} To conclude local asymptotical stability of the origin, it remains to show that $U_\ell$ is locally Lipschitz. Since $\sigma_\ell(V(\cdot))$ (resp. $W$) is locally Lipschitz, $U_\ell$ is locally Lipschitz in the region $W(\cdot)\leq\sigma_\ell(V(\cdot))$ (resp. $\sigma_\ell(V(\cdot))\leq W(\cdot)$). Since the hypotheses of Lemma \ref{thm:Dini derivative:stability by Dini derivative} (in Section \ref{sec:auxiliary results}) below are verified with $U(\cdot)=U_\ell(\cdot)$ and $E(\cdot)=E_\ell(\cdot)$, the origin is locally asymptotically stable for \eqref{eq:general system}. This concludes the proof of Proposition \ref{prop:asymetric:0 locally stable}.
\end{proof}

\subsection{Proof of Proposition \ref{prop:asymetric:neighborhood globally attractive}}\label{sec:proof of Proposition 2}

\begin{proof}
	This proof is analogous to the proof of Proposition \ref{prop:asymetric:0 locally stable} and divided into three parts. In the first one, it is obtained a function $\sigma_g$ that is class $\cK_\infty$ and $\cC^1$ with strictly positive derivative. This function is used in the second part, where a class $\cC^0$ proper and positive definite function $U_g$ is defined and its Dini derivative is studied. In the last part, it is used Lemma \ref{thm:Dini derivative:decreasing Dini derivative} to show that the set $\overline{\Omega_{M_g}(V)}\times\{0\}$ is globally asymptotically stable.

	{\itshape First part.} Consider the class $\cK$ functions $\delta$ and $\gamma_g$ from Assumptions \ref{hyp:asymmetric:z:ISS} and \ref{hyp:asymmetric:Lyapunov global}. The function $\gamma_{g}^{-1}$ is defined on $[0,b_g)$ and satisfies $\lim_{s\nearrow b_g}\gamma_{g}^{-1}(s)=\infty$. Assumption \ref{hyp:asymmetric:small gain assumption:global} implies that, $\forall s\in[M_g,b_g)$, $\delta(s)<\gamma_{g}^{-1}(s)$.
	Since $\gamma_g$ is of class $\cK$, from Lemma \ref{lem:construction of K_inf function} (in Section \ref{sec:auxiliary results}), there exists a class $\cK_\infty$ function $\tilde{\gamma}_g$ such that, $\forall s\in\Rs$,
	\begin{equation}\label{eq:proof:prop 0 locally stable:delta < t gamma_g}
		\delta(s)<\tilde{\gamma}_g(s)
	\end{equation}
	and, $\forall s\in[M_g,b_g)$,
	\begin{equation}\label{eq:proof:prop 0 locally stable:t gamma_g < gamma_g^-1}
		\tilde{\gamma}_g(s)<\gamma_g^{-1}(s).
	\end{equation}
	
	Since $\delta$ is of class $\cK$ and $\tilde{\gamma}_g$ is of class $\cK_\infty$ satisfying, $\forall s\in\Rs$, the inequality \eqref{eq:proof:prop 0 locally stable:delta < t gamma_g}, from Lemma \cite[Lemma A.1]{Jiangetal:1996}, there exists a function $\sigma_g$ that is of class $\cK_\infty$ and $\cC^1$ whose derivative is strictly positive and satisfies, $\forall s\in\Rs$,
\begin{equation}\label{eq:proof:prop 0 locally stable:delta < sigma_g < gamma_g^-1}
	\delta(s)<\sigma_g(s)<\tilde{\gamma}_g(s).
\end{equation}

	{\itshape Second part.} Let, $\forall (x,z)\in\R^n\times\R^m$, $U_g(x,z)\coloneqq\max\{\sigma_g(V(x)),W(z)\}$. Note that the function $U_g$ is proper positive definite. Pick $(x,z)\in\R^n\times\R^m$, one of three cases is possible: $\sigma_g(V({x}))<W({z})$, $W({z})<\sigma_g(V({x}))$ or $W({z})=\sigma_g(V({x}))$. The proof follows by showing that the Dini derivative of $U_g$ is negative. For each case, assume that $(x,z)\in(\R^n\setminus\Omega_{M_g}(V))\times\R^m$.

	\hspace{2pt}\underline{Case 1}. Suppose that $\sigma_g(V(x))<W(z)$. Analogously to the Case 1 of proof of Proposition \ref{prop:asymetric:0 locally stable}, $\sigma_g(V(x))<W(z)\Rightarrow D^+_{f,g}U_g(x,z)\leq-\lambda_z(z)$. This concludes Case 1.
	
	\hspace{2pt}\underline{Case 2}. Suppose that $W(z)<\sigma_g(V(x))$. This implies that
	$U_g(x,z)=\sigma_g(V(x))$ and $D^+_{f,g}U_g(x,z)=\sigma_g'(V(x))D_f^+V(x,z)$.
	From \eqref{eq:proof:prop 0 locally stable:delta < sigma_g < gamma_g^-1}, the following inequality
	\begin{equation}\label{eq:proof:prop neighborhood globally stable:W < sigma_g < t gamma_g}
		W(z)<\sigma_g(V(x))<\tilde{\gamma}_g(V(x))
	\end{equation}
	holds. At this point, two regions of $x$ will be analyzed: $b_g\leq V(x)$ and $M_g\leq V(x)< b_g$.
	
	\hspace{2pt}\underline{Case 2.a}. In the region where $b_g\leq V(x)$, Equation \eqref{eq:asymmetric:D_f^+V<0:global} together with \eqref{eq:asymmetric:b_g} yields $D_f^+V(x,z)\leq-\lambda_x(x)$.
	
	\hspace{2pt}\underline{Case 2.b}. In the region where $M_g\leq V(x)< b_g$, from \eqref{eq:proof:prop 0 locally stable:t gamma_g < gamma_g^-1} and \eqref{eq:proof:prop neighborhood globally stable:W < sigma_g < t gamma_g}, it yields
			\begin{equation}\label{eq:proof:prop 0 locally stable:W<delta<t gamma_g<gamma_g^-1}
					W(z)<\sigma_g(V(x))<\tilde{\gamma}_g(V(x))<\gamma^{-1}_g(V(x)).
			\end{equation}
			Equation \eqref{eq:asymmetric:D_f^+V<0:global} together with \eqref{eq:proof:prop 0 locally stable:W<delta<t gamma_g<gamma_g^-1} yields $D_f^+V(x,z)\leq-\lambda_x(x)$.
		
	Since, $\forall s\in\Rs$, $\sigma_g'(s)>0$, it follows that $D^+_{f,g}U_g(x,z)=\sigma_g'(V(x))D_f^+V(x,z)\leq-\sigma_g'(s)\lambda_x(x)$. This concludes Case 2.

	\hspace{2pt}\underline{Case 3}. Let $W(z)=\sigma_g(V(x))\coloneqq U_g^\ast(x,z)$. Analogously to the Case 3 of proof Proposition \ref{prop:asymetric:0 locally stable} and together with the analysis of Cases 1 and 2, the implication $W(z)=\sigma_g(V(x))\Rightarrow D^+_{f,g}U_g(x,z)\leq-\min\{\sigma_g'(s)\lambda_x(x),\lambda_z(z)\}$ holds, since $(x,z)\in(\R^n\setminus\Omega_{M_g}(V))\times\R^m$.

 Let $\tilde{M}_g=\min\{c\in\Rs:\overline{\Omega_{M_g}(V)}\times\{0\}\subset\overline{\Omega_{c}(U_g)}\ \text{and}\ \overline{\Omega_{c}(U_g)}\ \text{is connected}\}$. To sum up all the above cases, $\forall (x,z)\in(\R^n\times\R^m)\setminus\overline{\Omega_{\tilde{M}_g}(U_g)}$, 
	\begin{equation}\label{eq:proof:prop neighborhood globally stable:DU_g<0}
		\tilde{M}_g<U_g(x,z)\Rightarrow D^{+}_{f,g}U_g(x,z)\leq -E_g(x,z),
	\end{equation}
	where $E_g(\cdot,\cdot)=\min\{\sigma_g'(V(\cdot))\lambda_x(\cdot),\lambda_z(\cdot)\}$ is continuous and positive definite.
			
			{\itshape Third part.} Analogously to the third part of the proof of Proposition \ref{prop:asymetric:0 locally stable}, it follows that $U_g$ is locally Lipschitz. From Lemma \ref{thm:Dini derivative:equality of Dini time derivative and direction derivative} and \eqref{eq:proof:prop neighborhood globally stable:DU_g<0}, it follows that, $\forall (x,z)\in\R^n\times\R^m$ and $\forall t\in\Ras$, along solutions of \eqref{eq:general system},
			\begin{equation*}
				D^+U_g(X(t,x,z),Z(t,z,x))\hspace{-1.5pt}{\small=}\hspace{-1.5pt}D^+_{f,g}U_g(X(t,x,z),Z(t,z,x)).
			\end{equation*}
			
			Since solutions of \eqref{eq:general system} are absolutely continuous functions and, along solutions of \eqref{eq:general system}, $E_g$ is a continuous positive definite function, from Lemma \ref{thm:Dini derivative:decreasing Dini derivative}, $\forall (x,z)\in (\R^n\times\R^m)\setminus\overline{\Omega_{\tilde{M}_g}(U_g)}$ and $\forall t\in\Ras$, the function
			
			\begin{equation}\label{eq:proof:prop neighborhood globally stable:U_g<beta}
				t\mapsto U_g(X(t,x,z),Z(t,z,x))
			\end{equation}
			is strictly decreasing. Pick $(x,z)\in (\R^n\times\R^m)\setminus\overline{\Omega_{\tilde{M}_g}(U_g)}$, it will be proven that
			$$U^\infty_g\coloneqq\lim_{t\to\infty}U_g(X(t,x,z),Z(t,z,x))\leq\tilde{M}_g.$$
			To see the above suppose, by contradiction, that $U_g^\infty>\tilde{M}_g$. From the continuity of $U_g$, $\exists\varepsilon>0$ such that $U_g^\infty-\varepsilon>\tilde{M}_g$ and $U_g^\infty-\varepsilon\leq U_g(x,z)\leq U_g^\infty+\varepsilon$. Since $U_g$ is proper, the constant $\xi=\min\{E_g(x,z)>0:(x,z)\in U_g(x,z)\ \text{and}\ U_g^\infty-\varepsilon\leq U_g(x,z)\leq U_g^\infty+\varepsilon\}$ exists. Recalling the definition of $U_g$, $\exists T>0$ such that, $\forall t\geq T$, $U_g(X(t,x,z),Z(t,z,x))-U_g^\infty<\varepsilon$.  Moreover, from the definition of the constant $\xi$,
			\begin{equation*}
				\begin{split}
					U_g(X(t,x,z),Z(t,z,x))-U_g(X(T,x,z),Z(T,z,x))=\\
					\textstyle\int_{T}^t D^+U_g(X(s,x,z),Z(s,z,x))\,ds\leq -\xi\cdot(t-T).
				\end{split}
			\end{equation*}
			Then, 
			\begin{equation*}
				\begin{split}
					U^\infty_g=\lim_{t\to\infty}U_g(X(t,x,z),Z(t,z,x))\\
					=U_g(X(T,x,z),Z(T,z,x))\\
					+\lim_{t\to\infty}\textstyle\int_{T}^t D^+U_g(X(s,x,z),Z(s,z,x))\,ds\leq-\infty
				\end{split}
			\end{equation*}
		which contradicts the fact that $U_g$ is positive definite. Therefore, $U_g^\infty\leq \tilde{M}_\ell$. 
		
		In summary, the following facts hold for the function $U_g$: 1) $U_g$ is a proper positive definite function; 2) $U_g$ decreases along solutions of \eqref{eq:general system} having initial conditions in $(\R^n\times\R^m)\setminus\overline{\Omega_{\tilde{M}_g}(U_g)}$. From facts 1) and 2), the set $\overline{\Omega_{\tilde{M}_g}(U_g)}$ is globally asymptotically stable for \eqref{eq:general system}. This concludes the proof of Proposition \ref{prop:asymetric:neighborhood globally attractive}.
\end{proof}

\subsection{Proof of Theorem \ref{prop:asymetric:0 globally stable}}\label{sec:proof of Theorem 2}

\begin{proof}
	Under Assumption \ref{hyp:asymmetric:M_l<M_g}, $\exists M>0$ such that $M_g<M<M_\ell$. Under Assumptions \ref{hyp:asymmetric:z:ISS}, \ref{hyp:asymmetric:Lyapunov}, \ref{hyp:asymmetric:small gain assumption:local} and Proposition \ref{prop:asymetric:0 locally stable}, it follows that the origin is locally asymptotically stable. From the proof of Proposition \ref{prop:asymetric:0 locally stable}, there exists a proper positive definite function given, $\forall(x,z)\in\R^n\times\R^m$, by $U_\ell(x,z)=\max\{\sigma_\ell(V(x),W(z)\}$, where $\sigma_\ell$ is of class $\cK_\infty$ and $\cC^1$ satisfying \eqref{eq:proof:prop 0 locally stable:delta < sigma_l < gamma^-1}. Moreover, letting $\hat{M}_\ell\coloneqq\max\{c\in\Rs:c>M,\overline{\Omega_{c}(U_\ell)}\subset\overline{\Omega_{M_\ell}(V)}\times\{0\}\ \text{with}\ \overline{\Omega_{c}(U_\ell)} \ \text{connected}\}$ every solution starting in $\overline{\Omega_{\hat{M}_\ell}(U_\ell)}$ converges to the origin.

	Together with Assumptions \ref{hyp:asymmetric:z:ISS}, \ref{hyp:asymmetric:Lyapunov global}, \ref{hyp:asymmetric:small gain assumption:global} and the proof of Proposition \ref{prop:asymetric:neighborhood globally attractive}, it is possible to define, $\forall s\in\Ras$, a class $\cK_\infty$ function $\hat{\gamma}_g(s)=\min\{\tilde{\gamma}_g(s),\sigma_\ell(s)\}$ satisfying \eqref{eq:proof:prop 0 locally stable:delta < t gamma_g} and \eqref{eq:proof:prop 0 locally stable:t gamma_g < gamma_g^-1}. Then, it is obtained a class $\cK_\infty$ and $\cC^1$ function $\hat{\sigma}_g$ whose derivative is strictly positive and satisfies, $\forall s\in\Rs$,
	\begin{equation*}\tag{\ref{eq:proof:prop 0 locally stable:delta < sigma_g < gamma_g^-1}.new}
		\delta(s)<\hat{\sigma}_g(s)<\hat{\gamma}_g(s).
	\end{equation*}
	Defining a proper positive definite function given, $\forall(x,z)\in\R^n\times\R^m$, by $\hat{U}_g(x,z)=\max\{\hat{\sigma}_g(V(x)),W(z)\}$ and the constant $\hat{M}_g=\min\{c\in\Rs:c<M,\overline{\Omega_{M_g}(V)}\times\{0\}\subset\overline{\Omega_c(\hat{U}_g)} \ \text{with}\ \overline{\Omega_{c}(\hat{U}_g)} \ \text{connected}\}$, it follows from the proof of Proposition \ref{prop:asymetric:neighborhood globally attractive} that the set $\overline{\Omega_{\hat{M}_g}(\hat{U}_g)}$ is globally asymptotically stable.
	
	Since, $\forall s\in\Rs$, $\hat{\sigma}_g(s)<\sigma_\ell(s)$, it follows that, $\forall(x,z)\in(\R^n\times\R^m)\setminus\{(0,0)\}$, $\hat{U}_g(x,z)<U_\ell(x,z)$. This inequality implies that, $\forall c\in\Rs$, $\overline{\Omega_c(\hat{U}_g)}\subset\overline{\Omega_c(U_\ell)}$. Then, the following inclusion holds
	\begin{equation}\label{eq:inclusion}
		\overline{\Omega_{\hat{M}_g}(\hat{U}_g)}\subset\overline{\Omega_{M}(\hat{U}_g)}\subset\overline{\Omega_{M}(U_\ell)}\subset\overline{\Omega_{\hat{M}_\ell}(U_\ell)}.
	\end{equation}
	Thus, every solution of \eqref{eq:general system} starting in $(\R^n\times\R^m)\setminus\overline{\Omega_{\hat{M}_\ell}(U_\ell)}$ converges to $\overline{\Omega_{\hat{M}_g}(\hat{U}_g)}$, in finite time. Then, due to \eqref{eq:inclusion}, $\overline{\Omega_{\hat{M}_g}(\hat{U}_g)}\subset\overline{\Omega_{\hat{M}_\ell}(U_\ell)}$ holds, and thus solutions will converge to the origin, as $t\to\infty$. 
	
	From the above, combining the local asymptotical stability of the origin with its global attractivitty it is concluded that the origin is globally asymptotically stable for \eqref{eq:general system}.
\end{proof}

\section{Conclusion and perspectives}\label{sec:conclusion}

In this work, the authors shown that it is possible to make use of local and non-local input-to-state properties of an ISS system, in order to derive an ``optimal'' ISS gain. As a result of such approach, it is possible to apply the Small Gain Theorem in a less conservative way by deriving local and non-local small gain conditions to ensure the stability of an interconnected system. 

In a future work, the authors will generalize the above results for the case in which there exist four ISS gains: two for each subsystem. Moreover, the authors also intend to use the region-dependent gain condition to develop a methodology for the design of feedback laws under different gains constraints.

\section{Auxiliary results}\label{sec:auxiliary results}

\begin{lemma}\label{lem:construction of K_inf function}
	Let $\beta$ be a class $\cK$ function with
	\begin{equation}\label{eq:lemma:construction of K_inf function:lim beta}
		b=\lim_{s\to\infty}\beta(s).
	\end{equation}
	Let also $p,q$ be two constants and $\alpha$ be a class $\cK$ function such that, $0<p<q$ and, $\forall s\in[p,q]$,
	\begin{equation}\label{eq:lemma:construction of K_inf function:b circ a < s}
		\beta\circ\alpha(s)<s.
	\end{equation}
	Then, the class $\cK_\infty$ function $\tilde{\beta}$ given by
	\begin{equation}\label{eq:lemma:construction of K_inf function:t b}
		\tilde{\beta}(s)\hspace{-2pt}\coloneqq\hspace{-2pt} \left\{\begin{array}{l}
												\alpha(s)+\min\{s,K\},\text{if}\ p\neq0\ \text{and}\ s\in[0,p),\\
												\alpha(s)+\min\left\{s,\tfrac{\beta^{-1}(s)-\alpha(s)}{2}\right\},\text{if}\ q+\varepsilon<b\\ \hfill \text{and}\ s\in[p,q],\\
												A+B(s-q),\text{if}\ q+\varepsilon<b \ \text{and}\ s\in[q,q+\varepsilon),\\
												\alpha(s)+s,\text{if}\ q+\varepsilon\geq b\ \text{or}\ s\in[q+\varepsilon,\infty),
											\end{array}\right.
	\end{equation}
	 is such that, $\forall s\in\Rs$,
	\begin{equation}\label{eq:lemma:construction of K_inf function:a < t beta}
		\alpha(s)<\tilde{\beta}(s).
	\end{equation}
	Moreover, $\forall s\in[p,q]$, it also satisfies
	\begin{equation}\label{eq:lemma:construction of K_inf function:t beta < beta^-1}
		\tilde{\beta}(s)<\beta^{-1}(s).
	\end{equation}
\end{lemma}

Due to space constraints, the proof of Lemma \ref{lem:construction of K_inf function} is not provided in this paper.
	
	\begin{lemma}\label{thm:Dini derivative:stability by Dini derivative}\cite[Th\'eor\`eme 2.133]{Praly2011}
		Let $\bS\subset\R^k$ be a be a neighborhood of the origin. Let also the class $\cC^0$ function $h:\R^{k}\to\R^{k}$ and consider the system $\dot{y}=h(y)$. If there exist a positive definite and locally Lipschitz function $U:\bS\to\R$ and a positive definite function $E:\bS\to\R$ such that, $\forall y\in\bS$, $D^+_{h}U(y)\leq -E(y)$. Then, the origin is locally asymptotically stable for $\dot{y}=h(y)$.
	\end{lemma}
	
	\begin{lemma}\label{thm:Dini derivative:equality of Dini time derivative and direction derivative}\cite[Lemme 1.28]{Praly2011}
		Let the measurable and essentially bounded function $d:\R\to\R^p$ and the class $\cC^0$ function $h:\R^k\times\R^p\to\R^k$. If $U:\R^k\to\R$ is locally Lipschitz, then, for all maximal solutions $Y(t,y,d)$ of the system $\dot{y}=h(y,d(t))$ defined in the interval $(t_-,t_+)$, the function $t\mapsto U(Y(t,y,d))$, defined over $(t_-,t_+)$, is locally Lipschitz and, for almost every $ t\in(t_-,t_+)$,
		\begin{equation*}
			\tfrac{\partial U(Y)}{\partial t}(t,y,d)=D^+U(Y(t,y,d))=D^+_hU(Y(t,y,d)).
		\end{equation*}
		Moreover, if $d$ is continuous, the above equality holds, $\forall t\in(t_-,t_+)$.
	\end{lemma}
	
	
	\begin{lemma}\label{thm:Dini derivative:decreasing Dini derivative}
		Let $Y:\R\to\R^k$ be an absolutely continuous function, $U:\R^k\to\R$ be a locally Lipschitz proper positive definite function and $E:\R^k\to\R$ be a continuous positive definite function. Define, $\forall t\in\R$, $U(t)=U\circ Y(t)$ and $E(t)=E\circ Y(t)$. If, $\forall t\in\R$, $D^+U(t)\leq-E(t)$, then, $\forall t\in\R$, $U(t)$ is strictly decreasing.
	\end{lemma}
	
	Due to space constraints, the proof of Lemma \ref{thm:Dini derivative:decreasing Dini derivative} is not provided in this paper.
	
	\bibliographystyle{plain}
\bibliography{Library}
\end{document}